\documentclass[12pt,twoside]{article}
\usepackage[english]{babel}
\usepackage[latin1]{inputenc}
\usepackage{amsmath}
\usepackage{amssymb,amsfonts}
\usepackage{graphicx}
\usepackage{times,amssymb,amscd}
\newtheorem{thm}{Theorem}[section]

\newtheorem{lem}[thm]{Lemma}

\newtheorem{defn}[thm]{Definition}
\newtheorem{rem}[thm]{Remark}
\numberwithin{equation}{section}

\newcommand{\bA}{\mathbf{A}}

\newcommand{\bE}{\mathbf{E}}

\newcommand{\bH}{\mathbf{H}}

\newcommand{\bL}{\mathbf{L}}

\newcommand{\bR}{\mathbf{R}}
\newcommand{\bS}{\mathbf{S}}
\newcommand{\bV}{\mathbf{V}}

\newcommand{\be}{\mathbf{e}}

\newcommand{\br}{\mathbf{r}}
\newcommand{\bx}{\mathbf{x}}

\newcommand{\bT}{\mathbf{T}}

\newcommand{\bu}{\mathbf{u}}
\newcommand{\bv}{\mathbf{v}}
\newcommand{\bt}{\mathbf{t}}

\newcommand{\BV}{\boldsymbol{V}}

\newcommand{\Be}{\boldsymbol{e}}
\newcommand{\Bu}{\boldsymbol{u}}
\newcommand{\Bv}{\boldsymbol{v}}

\newcommand{\cP}{\mathcal{P}}
\newcommand{\cS}{\mathcal{S}}

\newcommand{\EUC}{\mathbf E^3}
\newcommand{\SPH}{\bS^3}
\newcommand{\HYP}{\bH^3}
\newcommand{\SXR}{\bS^2\!\times\!\bR}
\newcommand{\HXR}{\bH^2\!\times\!\bR}
\newcommand{\SLR}{\widetilde{\bS\bL_2\bR}}
\newcommand{\NIL}{\mathbf{Nil}}
\newcommand{\SOL}{\mathbf{Sol}}

\begin{document}
\pagestyle{myheadings}
\markboth{\centerline{Géza Csima and Jen\H o Szirmai}}
{Translation-like isoptic surfaces and angle sums $\dots$}
\title
{Translation-like isoptic surfaces and angle sums of translation triangles in $\NIL$ geometry
}

\author{Géza Csima and Jen\H o Szirmai \\
\normalsize Department of Geometry, Institute of Mathematics,\\
\normalsize Budapest University of Technology and Economics, \\
\normalsize M\"uegyetem rkp. 3., H-1111 Budapest, Hungary \\
\normalsize csimageza@gmail.com,~szirmai@math.bme.hu
\date{\normalsize{\today}}}

\maketitle
\begin{abstract}
After having investigated the geodesic and translation triangles and their angle 
sums in $\SOL$ and $\SLR$ geometries we consider the analogous problem 
in $\NIL$ space that
is one of the eight 3-dimensional Thurston geometries.

We analyze the interior angle sums of translation triangles in $\NIL$ geometry
and we provide a new approach to prove that it can be larger than or equal to $\pi$. 

Moreover, for the first time in non-constant curvature Thurston geometries 
we have developed a procedure for determining the equations of $\NIL$ isoptic
surfaces of translation-like segments and 
as a special case of this we examine the $\NIL$ translation-like Thales sphere, 
which we call {\it Thaloid}. 

In our work we will use the projective model of $\NIL$ described by E. Moln\'ar in \cite{M97}.
\end{abstract}

\section{Introduction} \label{section1}

In this paper we are interested in translation triangles and isoptic surfaces in $\NIL$ 
space that is one of the eight Thurston geometries (see \cite{S} 
and \cite{T}) derived by the Heisenberg matrix group \cite{M97, MSz}. 

In the Thurston spaces translation curves can be 
introduced in a natural way (see \cite{MoSzi10, Sz12}) by translations 
mapping each point to any point.
Consider a unit vector at the origin. Translations, postulated at the
beginning carry this vector to any point by its tangent mapping. If a curve $t\rightarrow (x(t),y(t),z(t))$ has just the translated
vector as tangent vector in each point, then the  curve is
called a {\it translation curve}. This assumption leads to a system of first order 
differential equations, thus translation
curves are simpler than geodesics and differ from them in $\NIL$, $\SLR$ and $\SOL$ 
geometries. Moreover, they play an important role and often 
seem to be more natural in these geometries, than their geodesic lines.  

In the remaining five Thurson geometries 
$\EUC,$ $\SPH,$ $\HYP,$ $\SXR$ and $\HXR,$ the translation and geodesic 
curves coincide with each other. 

Internal angle sum for triangles in $\SXR$ and $\HXR$ had been studied in \cite{Sz202}.

In \cite{CsSz16} we investigated the angle sum of translation and geodesic triangles in $\SLR$
geometry and proved that the possible sum of the interior angles in a 
translation triangle must be greater than or equal to $\pi.$ 
However, in geodesic triangles this sum can be less than, greater than or equal 
to $\pi.$

In \cite{Sz20} interior angle sum of translation triangles had been studied 
in $\SOL$ geometry to obtain that it must be greater than or equal 
to $\pi.$ To calculate this sum for geodesic curves needs further research.

In \cite{B} Brodaczewska showed, that the sum of the interior angles 
of translation triangles of $\NIL$ space is larger than or equal to $\pi,$ 
which is also the aim of this study. However, our approach below seems 
generally effective to all three geometries, 
where geodesic and translation curves differ. 


\begin{rem}
Of the Thurston geometries, those with constant curvature (Euclidean $\EUC$, hyperbolic $\HYP$, spherical $\SPH$)  have been extensively studied,  
but the other five geometries, $\HXR$, $\SXR$, $\NIL$, $\SLR$, $\SOL$ 
have been thoroughly studied only from a differential geometry and topological 
point of view. However, classical concepts highlighting the beauty and underlying 
structure of these geometries -- such as geodesic curves and spheres, translation curves an spheres, 
the lattices,  
the geodesic and translation triangles and their surfaces, their interior 
sum of angles, locus of points in the plane or in the space from a 
segment subtends a given angle (isoptic curves or surfaces) 
and similar statements to those 
known in constant curvature geometries -- can be formulated. These 
have not been the focus of attention (see \cite{Sz22,Sz22-3}). 

In this paper we consider some of these topics.
\end{rem}

In Section 2 we describe the projective model of $\NIL$
and we shall use its standard Riemannian metric obtained by pull back transform to the infinitesimal arc-length-square at the origin. We also recall the isometry group of $\NIL$ 
and give an overview about translation curves.

In Section 3 we study the $\NIL$ translation triangles and prove that the interior angle sum of a translation triangle in $\NIL$
geometry can be larger than, or equal to $\pi.$ 
We also determine when this internal angle sum is exactly $\pi.$

In Section 4 we introduce isoptic curves with the usual Euclidean planar definition, then we define the basic idea of spatial visibility, i.e. the concept of isoptic surfaces.
With the help of these, we can introduce the isoptic and orthoptic surfaces for any segment, which can then be defined for translation-like segments similarly in the $\NIL$ geometry.

In the last Section 5, we examine the translation-like isoptic surfaces of a given translation-like segment. We give the definition of these 
surfaces and a procedure by which we can determine the translation-like isoptic surface of any translation-like segment. 
This procedure can perhaps be further developed, to determine translation-like isoptic surfaces to any given $\NIL$ curve. 
We determine the implicit equation of these surfaces and visualize them. Thaloids are also analyzed as a special case. Similar 
investigations in this topic have only been carried out in spaces with constant curvature (see \cite{CsSz4,CsSz5,tizenegy,Odehnal,Cs-Sz14}).
\section{$\NIL$ geometry and its translation curves}
$\NIL$ geometry can be derived from the famous real matrix group $\mathbf{L(\mathbb{R})}$ discovered by Werner {Heisenberg}. The left (row-column) 
multiplication of Heisenberg matrices
     \begin{equation}
     \begin{gathered}
     \begin{pmatrix}
         1&x&z \\
         0&1&y \\
         0&0&1 \\
       \end{pmatrix}
       \begin{pmatrix}
         1&a&c \\
         0&1&b \\
         0&0&1 \\
       \end{pmatrix}
       =\begin{pmatrix}
         1&a+x&c+xb+z \\
         0&1&b+y \\
         0&0&1 \\
       \end{pmatrix}
      \end{gathered} \label{2.1}
     \end{equation}
defines "translations" $\mathbf{L}(\mathbb{R})= \{(x,y,z): x,~y,~z\in \mathbb{R} \}$ 
on the points of $\NIL= \{(a,b,c):a,~b,~c \in \mathbb{R}\}$. 
These translations are not commutative in general. The matrices $\mathbf{K}(z) \vartriangleleft \mathbf{L}$ of the form
     \begin{equation}
     \begin{gathered}
       \mathbf{K}(z) \ni
       \begin{pmatrix}
         1&0&z \\
         0&1&0 \\
         0&0&1 \\
       \end{pmatrix}
       \mapsto (0,0,z)  
      \end{gathered}\label{2.2}
     \end{equation} 
constitute the one parametric centre, i.e. each of its elements commutes with all elements of $\mathbf{L}$. 
The elements of $\mathbf{K}$ are called {\it fibre translations}. $\NIL$ geometry of the Heisenberg group can be projectively 
(affinely) interpreted by "right translations" 
on points as the matrix formula 
     \begin{equation}
     \begin{gathered}
       (1;a,b,c) \to (1;a,b,c)
       \begin{pmatrix}
         1&x&y&z \\
         0&1&0&0 \\
         0&0&1&x \\
         0&0&0&1 \\
       \end{pmatrix}
       =(1;x+a,y+b,z+bx+c) 
      \end{gathered} \label{2.3}
     \end{equation} 
shows, according to \ref{2.1}. Here we consider $\mathbf{L}$ as projective collineation group with right actions in homogeneous coordinates.
We will use the Cartesian homogeneous coordinate simplex $E_0(\be_0)$,$E_1^{\infty}(\be_1)$,$E_2^{\infty}(\be_2)$,
$E_3^{\infty}(\be_3), \ (\{\be_i\}\subset \bV^4$ \ $\text{with the unit point}$ $E(\be = \be_0 + \be_1 + \be_2 + \be_3 ))$ 
which is distinguished by an origin $E_0$ and by the ideal points of coordinate axes, respectively. 
Moreover, $\mathbf{y}=c\bx$ with $0<c\in \mathbb{R}$ (or $c\in\mathbb{R}\setminus\{0\})$
defines a point $(\bx)=(\mathbf{y})$ of the projective 3-sphere $\mathcal{P} \mathcal{S}^3$ (or that of the projective space $\cP^3$ where opposite rays
$(\bx)$ and $(-\bx)$ are identified). 
The dual system $\{(\Be^i)\}, \ (\{\Be^i\}\subset \BV_4)$, with $\be_i\Be^j=\delta_i^j$ (the Kronecker symbol), describes the simplex planes, especially the plane at infinity 
$(\Be^0)=E_1^{\infty}E_2^{\infty}E_3^{\infty}$, and generally, $\Bv=\Bu\frac{1}{c}$ defines a plane $(\Bu)=(\Bv)$ of $\cP \cS^3$
(or that of $\cP^3$). Thus $0=\bx\Bu=\mathbf{y}\Bv$ defines the incidence of point $(\bx)=(\mathbf{y})$ and plane
$(\Bu)=(\Bv)$, as $(\bx) \text{I} (\Bu)$ also denotes it. Thus {\bf Nil} can be visualized in the affine 3-space $\bA^3$
(so in $\bE^3$) as well \cite{MSz06}.

In this context E. Moln\'ar \cite{M97} has derived the well-known infinitesimal arc-length square invariant under translations $\bL$ at any point of $\NIL$ as follows
\begin{equation}
   \begin{gathered}
      (dx)^2+(dy)^2+(-xdy+dz)^2=\\
      (dx)^2+(1+x^2)(dy)^2-2x(dy)(dz)+(dz)^2=:(ds)^2
       \end{gathered} \label{2.4}
     \end{equation}
The translation group $\mathbf{L}$ defined by formula \ref{2.3} can be extended to a larger group $\mathbf{G}$ of collineations,
preserving the fibres, that will be equivalent to the (orientation preserving) isometry group of $\NIL$. 

In \cite{M06} E.~Moln\'ar has shown that 
a rotation through angle $\omega$
about the $z$-axis at the origin, as isometry of $\NIL$, keeping invariant the Riemann
metric everywhere, will be a quadratic mapping in $x,y$ to $z$-image $\overline{z}$ as follows:
     \begin{equation}
     \begin{gathered}
       \mathcal{M}=\br(O,\omega):(1;x,y,z) \to (1;\overline{x},\overline{y},\overline{z}); \\ 
       \overline{x}=x\cos{\omega}-y\sin{\omega}, \ \ \overline{y}=x\sin{\omega}+y\cos{\omega}, \\
       \overline{z}=z-\frac{1}{2}xy+\frac{1}{4}(x^2-y^2)\sin{2\omega}+\frac{1}{2}xy\cos{2\omega}.
      \end{gathered} \label{2.5}
     \end{equation}
This rotation formula $\mathcal{M}$, however, is conjugate by the quadratic mapping $\alpha$ to the linear rotation $\Omega$ as follows
     \begin{equation}
     \begin{gathered}
       \alpha^{-1}: \ \ (1;x,y,z) \stackrel{\alpha^{-1}}{\longrightarrow} (1; x',y',z')=(1;x,y,z-\frac{1}{2}xy) \ \ \text{to} \\
       \Omega: \ \ (1;x',y',z') \stackrel{\Omega}{\longrightarrow} (1;x",y",z")=(1;x',y',z')
       \begin{pmatrix}
         1&0&0&0 \\
         0&\cos{\omega}&\sin{\omega}&0 \\
         0&-\sin{\omega}&\cos{\omega}&0 \\
         0&0&0&1 \\
       \end{pmatrix}, \\
       \text{with} \ \ \alpha: (1;x",y",z") \stackrel{\alpha}{\longrightarrow}  (1; \overline{x}, \overline{y},\overline{z})=(1; x",y",z"+\frac{1}{2}x"y").
      \end{gathered} \label{2.6}
     \end{equation}
This quadratic conjugacy modifies the $\NIL$ translations in \ref{2.3}, as well. Now a translation with $(X,Y,Z)$ in \ref{2.3} instead of $(x,y,z)$ will be changed 
by the above conjugacy to the translation 
     \begin{equation}
     \begin{gathered}
        (1;x,y,z) \longrightarrow (1; \overline{x}, \overline{y},\overline{z})=(1; x,y,z)
        \begin{pmatrix}
         1&X&Y&Z-\frac{1}{2}XY \\
         0&1&0&-\frac{1}{2}Y \\
         0&0&1&\frac{1}{2}X \\
         0&0&0&1 \\
       \end{pmatrix}, \\
             \end{gathered} \label{2.7}
     \end{equation}
that is again an affine collineation.
\subsection{Translation curve and sphere}
We consider a $\NIL$ curve $(1,x(t), y(t), z(t) )$ with a given starting tangent vector at the origin $O=E_0=(1,0,0,0)$
\begin{equation}
   \begin{gathered}
      u=\dot{x}(0),\ v=\dot{y}(0), \ w=\dot{z}(0).
       \end{gathered} \label{2.8}
     \end{equation}
For a translation curve let its tangent vector at the point $(1,x(t), y(t), z(t) )$ be defined by the matrix \ref{2.3} 
with the following equation:
\begin{equation}
     \begin{gathered}
     (0,u,v,w)
     \begin{pmatrix}
         1&x(t)&y(t)&z(t) \\
         0&1&0&0 \\
         0&0&1&x(t) \\
         0&0&0&1 \\
       \end{pmatrix}
       =(0,\dot{x}(t),\dot{y}(t),\dot{z}(t)).
       \end{gathered} \label{2.9}
     \end{equation}
Thus, the {\it translation curves} in $\NIL$ geometry (see  \cite{MoSzi10}, \cite{MSz06} \cite{MSzV}) are defined by the above first order differential equation system 
$\dot{x}(t)=u, \ \dot{y}(t)=v,  \ \dot{z}(t)=v \cdot x(t)+w,$ whose solution is the following: 
\begin{equation}
   \begin{gathered}
       x(t)=u t, \ y(t)=v t,  \ z(t)=\frac{1}{2}uvt^2+wt.
       \end{gathered} \label{2.10}
\end{equation}
We assume that the starting point of a translation curve is the origin, because we can transform a curve into an 
arbitrary starting point by translation \ref{2.3}, moreover, unit initial velocity translation can be assumed by "geographic" parameters
$\phi$ and $\theta$:
\begin{equation}
\begin{gathered}
        x(0)=y(0)=z(0)=0; \\ \ u=\dot{x}(0)=\cos{\theta} \cos{\phi}, \ \ v=\dot{y}(0)=\cos{\theta} \sin{\phi}, \ \ w=\dot{z}(0)=\sin{\theta}; \\ 
        - \pi \leq \phi \leq \pi, \ -\frac{\pi}{2} \leq \theta \leq \frac{\pi}{2}. \label{2.11}
\end{gathered}
\end{equation}
\begin{defn}
The translation distance $d^t(P_1,P_2)$ between the points $P_1$ and $P_2$ is defined by the arc length of the above translation curve 
from $P_1$ to $P_2$.
\end{defn}

\begin{defn} The sphere of radius $r >0$ with centre at the origin, (denoted by $S^t_O(r)$), with the usual longitude and altitude parameters 
$\phi$ and $\theta$, respectively by \ref{2.11}, is specified by the following equations:
\begin{equation}
\begin{gathered}
        S_O^t(r): \left\{ \begin{array}{ll} 
        x(\phi,\theta)=r \cos{\theta} \cos{\phi}, \\
        y(\phi,\theta)=r \cos{\theta} \sin{\phi}, \\
        z(\phi,\theta)=\frac{r^2}{2} \cos^2{\theta} \cos{\phi} \sin{\phi}+ r \sin{\theta}.
        \end{array} \right.
        \label{2.12}
\end{gathered}
\end{equation}
\end{defn}
\begin{defn}
 The body of the translation sphere of centre $O$ and of radius $r$ in the $\NIL$ space is called translation ball, denoted by $B^t_{O}(r)$,
 i.e. $Q \in B^t_{O}(r)$ iff $0 \leq d^t(O,Q) \leq r$.
\end{defn}

The parametrization in \ref{2.12} allows us, to create the implicit equation of $B^t_{O}(r)$:
 \begin{equation}
\begin{gathered}
        x^2+y^2+\left(z-\dfrac{xy}{2}\right)^2=r^2
						\label{2.13}
\end{gathered}
\end{equation}

\section{Translation triangles}
We consider $3$ points $A_1$, $A_2$, $A_3$ in the projective model of $\NIL$ space.
The {\it translation segments} $a_k$ connecting the points $A_i$ and $A_j$
$(i<j,~i,j,k \in \{1,2,3\}, k \ne i,j$) are called sides of the {\it translation triangle} with vertices $A_1$, $A_2$, $A_3$.
\begin{figure}[ht]
\centering
\includegraphics[width=14cm]{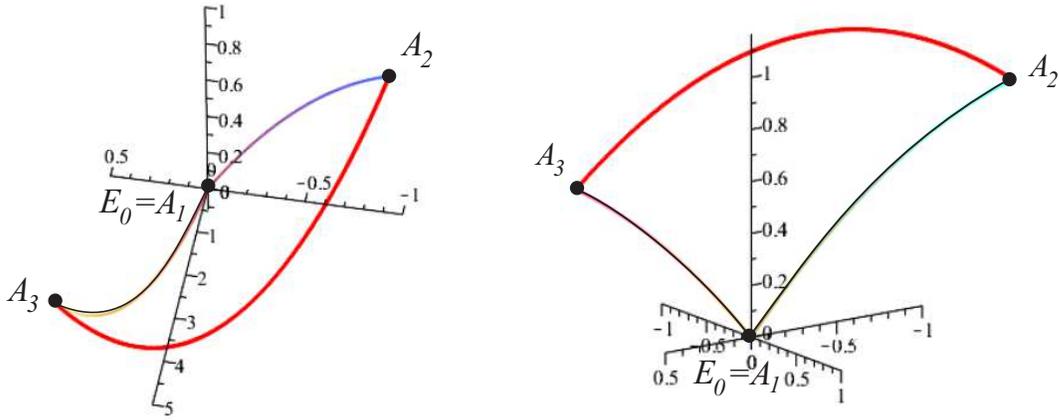}
\caption{Translation triangles with vertices $A_1=(1,0,0,0)$, $A_2=(1,-1,1,1)$, $A_3=(1,1/2,5,1/2)$ (left) and
with vertices $A_1=(1,0,0,0)$, $A_2=(1,-1,1,1)$, $A_3=(1,1/2,-1,1/2)$ (right).}
\label{fig:Fig1}
\end{figure}
In Riemannian geometries the metric tensor (or infinitesimal arc-lenght square (see \ref{2.4}) is used to define the angle $\theta$ between two curves.
If their tangent vectors in their common point are $\bu$ and $\bv$ and $g_{ij}$ are the components of the metric tensor then
\begin{equation}
\cos(\theta)=\frac{u^i g_{ij} v^j}{\sqrt{u^i g_{ij} u^j~ v^i g_{ij} v^j}} \label{3.1}
\end{equation}
It is clear by the above definition of the angles and by the infinitesimal arc-lenght square \ref{2.4}, that
the angles are the same as the Euclidean ones at the origin of
the projective model of $\NIL$ geometry.

Considering a translation triangle $A_1A_2A_3$ we can assume by the homogeneity of the $\NIL$ geometry that one of its vertex 
coincide with the origin $A_1=E_0=(1,0,0,0)$ and the other two vertices are $A_2(1,x^2,y^2,z^2)$ and $A_3(1,x^3,y^3,z^3)$. 

We will consider the {\it interior angles} of translation triangles that are denoted at the vertex $A_i$ by $\omega_i$ $(i\in\{1,2,3\})$.
We note here that the angle of two intersecting translation curves depends on the orientation of their tangent vectors. 

\begin{figure}[ht]
\centering
\includegraphics[width=13cm]{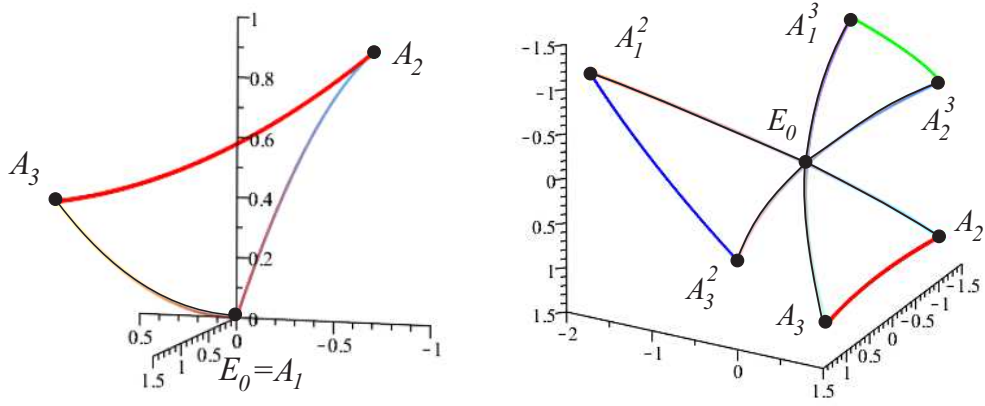}
\caption{Translation triangle with vertices $A_1=(1,0,0,0)$, $A_2=(1,-1,1,1)$, $A_3=(1,1/2,3/2,1/2)$ and its translated copies $A_1^2A_3^2E_0$ and $A_1^3A_2^3E_0$.}
\label{fig:Fig2}
\end{figure}

{\it In order to determine the interior angles of a translation triangle $A_1A_2A_3$ 
and its interior angle sum $\sum_{i=1}^3(\omega_i)$,
we define {translations} $\bT_{A_i}$, $(i\in \{2,3\})$ as elements of the isometry group of $\NIL$, that
maps the origin $E_0$ onto $A_i$} (see Fig. \ref{fig:Fig2}).

E.g. the isometry $\bT_{A_2}$ and its inverse (up to a positive determinant factor) can be given by:
\begin{equation}
\bT_{A_2}=
\begin{pmatrix}
1 & x_2 & y_2 & z_2 \\
0 & 1 & 0 & 0 \\
0 & 0 & 1 & x_2 \\
0 & 0 & 0 & 1
\end{pmatrix} , ~ ~ ~
\bT_{A_2}^{-1}=
\begin{pmatrix}
1 & -x_2 & -y_2 & x_2 y_2 - z_2 \\
0 & 1 & 0 & 0 \\
0 & 0 & 1 & -x_2 \\
0 & 0 & 0 & 1 
\end{pmatrix}, \notag
\end{equation}
\begin{equation}
\bT_{A_3}=
\begin{pmatrix}
1 & x_3 & y_3 & z_3 \\
0 & 1 & 0 & 0 \\
0 & 0 & 1 & x_3 \\
0 & 0 & 0 & 1
\end{pmatrix} , ~ ~ ~
\bT_{A_3}^{-1}=
\begin{pmatrix}
1 & -x_3 & -y_3 & x_3 y_3 -z_3 \\
0 & 1 & 0 & 0 \\
0 & 0 & 1 & -x_3 \\
0 & 0 & 0 & 1 , 
\end{pmatrix},\label{3.2}
\end{equation}
and the images $\bT^{-1}_{A_2}(A_i)$ of the vertices $A_i$ $(i \in \{1,2,3\})$ are the following (see also Fig. \ref{fig:Fig2}):
\begin{equation}
\begin{gathered}
T^{-1}_{A_{2}}(A_{1})=A^{2}_{1}=(1,-x_2,-y_2,x_2 y_2-z_2); \text{ } T^{-1}_{A_{2}}(A_2)=A^{2}_{2}=E_0=(1,0,0,0);\\
T^{-1}_{A_{2}}(A_3)=A^{2}_{3}=(1,-x_2+x_3,-y_2+y_3,x_2 y_2-x_2 y_3-z_2+z_3);\\
T^{-1}_{A_{3}}(A_{1})=A_{1}^{3}=(1,-x_3,-y_3,x_3 y_3 -z_3); \text{ } T^{-1}_{A_{3}}(A_{3})=A_{3}^{3}=E_0=(1,0,0,0); \\
T^{-1}_{A_{3}}(A_2)=A_{2}^{3}=(1,x_2-x_3,y_2-y_3,-x_3 y_2+x_3 y_3 +z_2 -z_3). \label{3.3}
\end{gathered}
\end{equation}

Our aim is to determine angle sum $\sum_{i=1}^3(\omega_i)$ of the interior 
angles of translation triangles $A_1A_2A_3$ (see Fig. \ref{fig:Fig1} and Fig. \ref{fig:Fig2}).
We have seen that $\omega_1$ and the angle of translation curves 
with common point at the origin $E_0$ is the same as the 
Euclidean one therefore can be determined by usual Euclidean sense.

The translations $\bT_{A_i}$ $(i=2,3)$ are isometries in $\NIL$ geometry thus
$\omega_i$ is equal to the angle $(t(A_i^i, A_1^i)t(A_i^i, A_j^i))\angle$ $(i,j=2,3$, $i \ne j)$ (see Fig. \ref{fig:Fig2})
where $t(A_i^i, A_1^i)$, $t(A_i^i, A_j^i)$ are oriented translation curves $(E_0=A_2^2=A_3^3)$ and
$\omega_1$ is equal to the angle $(t(E_0, A_2)t(E_0, A_3)) \angle$ 
where $t(E_0, A_2)$, $t(E_0, A_3)$ are also oriented translation curves.

We denote the oriented unit tangent vectors of the oriented geodesic curves $t(E_0, A_i^j)$ with $\mathbf{t}_i^j$ where
$(i,j)\in\{(1,3),(1,2),(2,3),(3,2),(3,0),(2,0)\}$ and $A_3^0=A_3$, $A_2^0=A_2$.
The Euclidean coordinates of $\mathbf{t}_i^j$ are :
\begin{equation}
t_{i}^{j}=(\cos{\theta_{i}^{j}} \cos{\phi_{i}^{j}},\cos{\theta_{i}^{j}} \sin{\phi_{i}^{j}},\sin{\theta_{i}^{j}}). \label{3.4}
\end{equation}
In order to obtain the angle of two translation curves $t_{E_0A_i^j}$ and $t_{E_0A_k^l}$ ($(i,j)\ne(k,l)$; $(i,j),(k,l)\in\{(1,3),(1,2),(2,3),(3,2),(3,0),(2,0)\})$ intersected at the origin $E_0$ we need to determine their tangent vectors $\bt_s^r$ 
$((s,r) \in \{(1,3),(1,2),$ $(2,3),(3,2),(3,0),(2,0)\})$
(see \ref{3.4}) at their starting point $E_0$.
From \ref{3.4} follows that a tangent vector at the origin is given 
by the parameters $\phi$ and $\theta$ of the corresponding translation 
curve (see \ref{2.11}) that 
can be determined from the homogeneous coordinates 
of the endpoint of the translation curve.

It can be assumed by the homogeneity of $\NIL$ that the starting point of a 
given translation curve segment is $E_0=P_1=(1,0,0,0)$ and 
the other endpoint will be given by its homogeneous coordinates $P_2=(1,a,b,c)$. 
We consider the translation curve segment $t_{P_1P_2}$ and determine its
parameters $(\phi,\theta,r)$ expressed by the real coordinates $a$, $b$, $c$ of $P_2$. 
We obtain directly by equation system \ref{2.11} the following:
\begin{lem}
\begin{enumerate}
\item Let $(1,a,b,c)$ $(a,b \in \mathbb{R} \setminus \{0\},~c\in \mathbb{R})$ be the homogeneous 
coordinates of the point $P \in \NIL$. The parameters of the
corresponding translation curve $t_{E_0P}$ are the following
\begin{equation}
\begin{gathered}
\phi=\mathrm{arccot}\Big(\frac{a}{b}\Big),~\text{or}~ \phi=\mathrm{arccot}\Big(\frac{a}{b}\Big)-\pi, \\
\theta=\mathrm{arctan}\Big( \frac{c-\frac{ab}{2}}{\sqrt{a^{2}+b^{2}}}\Big),~
r=\Big|\frac{c-\frac{ab}{2}}{\sin{\theta}}\Big|.
\end{gathered} \label{3.5}
\end{equation}
\item Let $(1,a,0,c)$ $(a,c \in \mathbb{R} \setminus \{0\})$ be the homogeneous 
coordinates of the point $P \in \NIL$. The parameters of the
corresponding translation curve $t_{E_0P}$ are the following
\begin{equation}
\begin{gathered}
\phi=\pi \cdot n ,~ (n\in\{0,1\}),~ 
\theta= \mathrm{arctan}\Big(\frac{c}{a}\Big),~
r=\Big|\frac{a}{\cos{\theta}}\Big|.
\end{gathered} \label{3.6}
\end{equation}
\item Let $(1,a,0,0)$ $(a \in \mathbb{R}\setminus \{0\})$ 
be the homogeneous coordinates of the point $P \in \NIL$. The parameters of the 
corresponding translation curve $t_{E_0P}$ are the following
\begin{equation}
\begin{gathered}
\phi=\pi \cdot n ,~ (n\in\{0,1\}),~
\theta=\pi \cdot n ,~ (n\in\{0,1\}),~
r=|a|.
\end{gathered} \label{3.7}
\end{equation}
\item Let $(1,0,b,0)$ $(b \in \mathbb{R}\setminus \{0\})$ 
be the homogeneous coordinates of the point $P \in \NIL$. The parameters of the 
corresponding translation curve $t_{E_0P}$ are the following
\begin{equation}
\begin{gathered}
\phi=\pm \frac{\pi}{2},~
\theta=\pi \cdot n ,~ (n\in\{0,1\}),~
r=|b|.
\end{gathered} \label{3.8}
\end{equation}
\item Let $(1,0,0,c)$ $(c \in \mathbb{R}\setminus \{0\})$ 
be the homogeneous coordinates of the point $P \in \NIL$. The parameters of the 
corresponding translation curve $t_{E_0P}$ are the following
\begin{equation}
\begin{gathered}
\theta=\pm \frac{\pi}{2},~
r=|c|.~ ~ \square
\end{gathered} \label{3.9}
\end{equation}
\end{enumerate}
\label{lem1}
\end{lem}
Applying the above lemma we obtain the following 
\begin{thm}
The sum of the interior angles of a translation triangle is greater than or equal to $\pi$.
\end{thm}
\textbf{Proof:} The translations $\bT_{A_2}^{-1}$ and $\bT_{A_3}^{-1}$ are isometries 
in $\NIL$ geometry thus $\omega_2$ is equal to the angle 
$((A_2^2 A_1^2), (A_2^2 A_3^2)) \angle$ (see Fig. \ref{fig:Fig2}) 
of the oriented translation segments $t_{A_2^2 A_1^2}$, 
$t_{A_2^2A_3^2}$ and $\omega_3$ is equal to the angle 
$((A_3^3 A_1^3),(A_3^3 A_2^3)) \angle$ 
of the oriented translation segments 
$t_{A_3^3 A_1^3}$ and $t_{A_3^3 A_2^3}$ $(E_0=A_2^2=A_3^3$).  

Substituting the coordinates of the points 
$A_i^j$ (see \ref{3.3} and \ref{3.4}) $((i,j) \in \{(1,3),(1,2),$ $(2,3),(3,2),(3,0),(2,0)\})$ to the appropriate equations of Lemma \ref{lem1}, 
it is easy to see that 
\begin{equation}
\begin{gathered}
\theta_2^0=-\theta_1^2,~\phi_2^0-\phi_1^2=\pm \pi \Rightarrow \bt_2^0=-\bt_1^2,\\ 
\theta_3^0=-\theta_1^3,~\phi_3^0-\phi_1^3=\pm \pi \Rightarrow \bt_3^0=-\bt_1^3,\\
\theta_3^2=-\theta_2^3,~\phi_3^2-\phi_2^3=\pm \pi \Rightarrow \bt_3^2=-\bt_2^3.
\end{gathered} \label{3.10}
\end{equation}
\begin{figure}[ht]
\centering
\includegraphics[width=8cm]{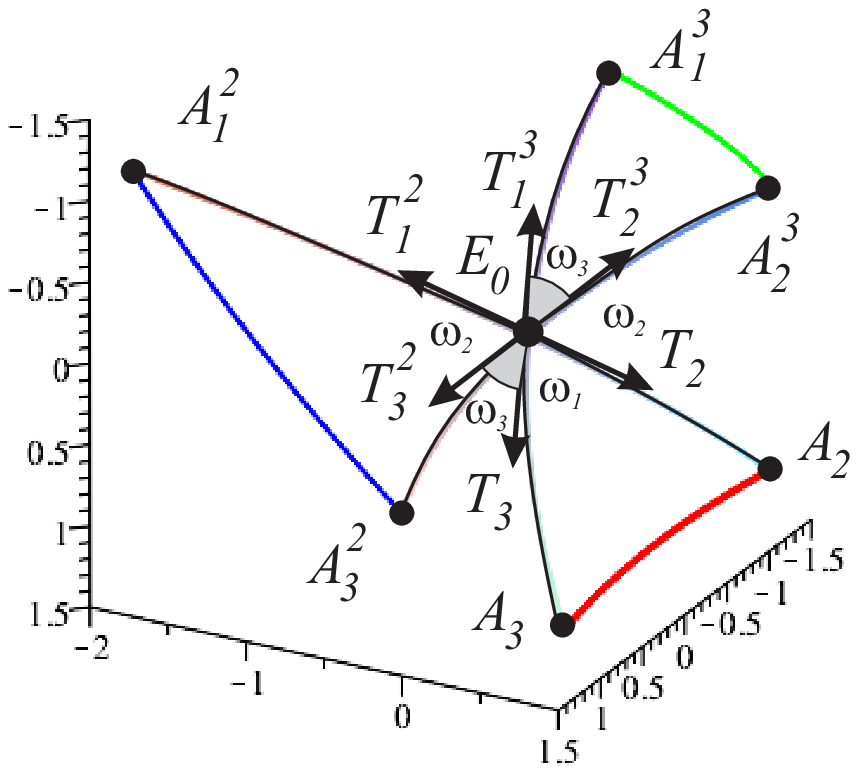}
\caption{Translation triangle with vertices $A_1=(1,0,0,0)$, $A_2=(1,-1,1,1)$, $A_3=(1,1/2,3/2,1/2)$, 
its translated copies $A_1^2A_3^2E_0$, $A_1^3A_2^3E_0$ with angles $\omega_i$ $(i\in\{1,2,3\}$.}
\label{Fig3}
\end{figure}
\begin{figure}
\centering
\includegraphics[width=13cm]{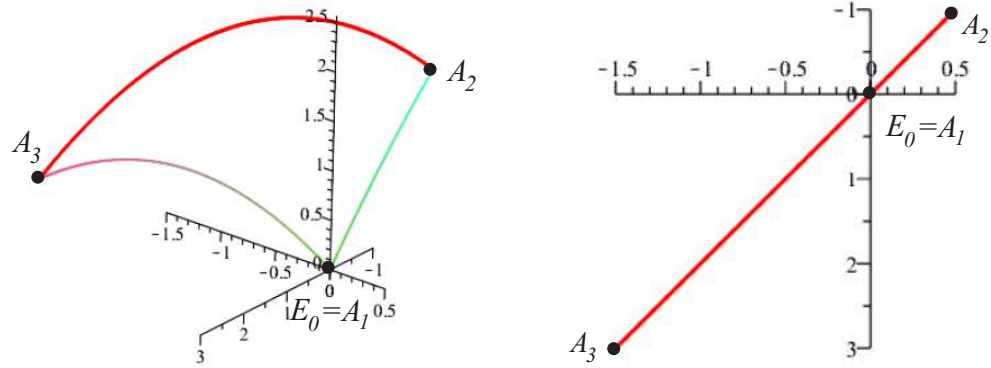}
\caption{Translation triangle with vertices $A_1=(1,0,0,0)$, $A_2=(1,-1,1/2,2)$, $A_3=(1,3,-3/2,1)$. 
The translation curve segments $t_{A_1A_2}$, $t_{A_2A_3}$, $t_{A_3A_1}$ lie
on a plane orthogonal to the $[x,y]$ base plane. The interior angle sum of this translation triangle is $\sum_{i=1}^3(\omega_i)=\pi$.}
\label{Fig4}
\end{figure}
The endpoints $T_i^j$ of the position vectors $\bt_i^j=\overrightarrow{E_0T_i^j}$  
lie on the unit sphere centered at the origin. The measure of angle $\omega_i$ $(i\in \{1,2,3\})$ of the vectors $\bt_i^j$ and $\bt_r^s$ is equal to the spherical 
distance of the corresponding points $T_i^j$ and $T_r^s$ on the unit sphere (see Fig.\ref{Fig3}). Moreover, a direct consequence of equations (\ref{3.9}) that each point pair 
($T_2$, $T_1^2$), $(T_3$,$T_1^3$), ($T_2^3$,$T_3^2$) 
contains antipodal points related to the unit sphere with centre $E_0$.

Due to the antipodality $\omega_1=T_2E_0T_3 \angle =T_1^2E_0T_1^3 \angle$, therefore their corresponding spherical 
distances are equal, as well (see Fig. \ref{Fig3}). 
Now, the sum of the interior angles $\sum_{i=1}^3(\omega_i)$ can be considered as three consecutive spherical arcs $(T_3^2 T_1^2)$, $(T_1^2 T_1^3)$, 
$T_1^3 T_2^3)$. 
Since the triangle inequality holds on the sphere, the sum of these arc lengths is greater or equal to the half 
of the circumference of the main circle on the unit sphere {i.e.} $\pi$. $\square$
\medbreak 
The following lemma is an immediate consequence of the above proof:
\begin{lem}
The angle sum $\sum_{i=1}^3(\omega_i)$ of a $\NIL$ translation
triangle $A_1A_2A_3$ is $\pi$ if and only if the points $T_i^j$ 
$((i,j) \in \{(1,3),(1,2),$ $(2,3),(3,2),(3,0),(2,0)\})$ lie 
in an Euclidean plane (Fig. \ref{Fig4}).
\label{lem2}
\end{lem}

Now, we distinguish the cases, when the internal angle sum is exactly $\pi.$

\begin{lem}
If the vertices of a translation triangle $A_1A_2A_3$ lie in a plane perpendicular 
to the base plane (coordinate plane $[x,y]$) of 
the model of $\NIL$ geometry 
then the interior angle sum $\sum_{i=1}^3(\omega_i)=\pi$. 
\label{lem3}
\end{lem}
{\bf{Proof}}: We get from the equation system (\ref{2.10}) of the translation curves that 
the points of a translation curve $t_{E_0P}$ ($P\in \NIL$) lie in an Euclidean plane 
that is perpendicular to $[x,y]$ base plane, therefore, its tangent 
line also lies in this plane.

Moreover, a direct consequence of formulas (\ref{2.10}) and (\ref{3.3}) than if a 
translation triangle $A_1A_2A_3$ lies in this to base plane perpendicular 
$\alpha$ plane then its translated image by a translation 
lies also in a to the base plane perpendicular plane $\alpha'$ and each to the base plane orthogonal $\alpha'$ 
can be derived as a translated copy of $\alpha$. Thus,
applying the Lemma \ref{lem2} we proved this lemma. \ \   $\square$
\medbreak
We can determine the interior angle sum of arbitrary translation triangle.
In the following table we summarize some numerical data of interior angles of given translation triangles:
\medbreak
\medbreak
\centerline{\vbox{
\halign{\strut\vrule\quad \hfil $#$ \hfil\quad\vrule
&\quad \hfil $#$ \hfil\quad\vrule &\quad \hfil $#$ \hfil\quad\vrule &\quad \hfil $#$ \hfil\quad\vrule &\quad \hfil $#$ \hfil\quad\vrule
\cr
\noalign{\hrule}
\multispan5{\strut\vrule\hfill\bf Table 2: ~ $A_2(1,-1,1,1),$ ~ $A_3(1,1/2,y^3,1/2)$ \hfill\vrule}%
\cr
\noalign{\hrule}
\noalign{\vskip2pt}
\noalign{\hrule}
y^3 & \omega_1 & \omega_2 & \omega_3  & \sum_{i=1}^3(\omega_i)  \cr
\noalign{\hrule}
-10 & 1.85298 & 0.75454 & 0.58205 & 3.18956 \cr
\noalign{\hrule}
-2 & 1.78411 &0.52781 & 0.83577 & 3.14770 \cr
\noalign{\hrule}
-1 & 1.70632 & 0.44929 & 0.98637 &  3.14198 \cr
\noalign{\hrule}
1/10 & 1.35152 & 0.46598 & 1.32927  & 3.14677 \cr
\noalign{\hrule}
3/4 & 1.19668 & 0.68254 & 1.31811 & 3.19733 \cr
\noalign{\hrule}
3/2 & 1.19912 & 1.08556 & 0.97181  & 3.25650 \cr
\noalign{\hrule}
5 & 1.24271 & 1.94607 & 0.36983  & 3.55861 \cr
\noalign{\hrule}
10 & 1.25686 & 2.12780 & 0.40324 &  3.78790 \cr
\noalign{\hrule}
}}}
\medbreak

\section{Introduction to isoptic curves}

It is well known that in the Euclidean plane the locus of points from a 
segment subtends a given angle $\alpha$ $(0<\alpha<\pi)$ is the union of 
two arcs except for the endpoints with the segment as common chord. If this $\alpha$ is equal to $\frac{\pi}{2}$ then we get the Thales circle. 
Replacing the segment to another general curve, we obtain the Euclidean definition of 
isoptic curve:
\begin{defn}[\cite{yates}]The locus of the intersection of tangents
to a curve meeting at a constant angle $\alpha$ $(0<\alpha <\pi)$ is
the $\alpha$ -- isoptic of the given curve. The isoptic 
curve with right angle called \textit{orthoptic curve}.
\label{defiso}
\end{defn}
\begin{rem}Sometimes we consider the $\alpha$ -- and $\pi-\alpha$ -- isoptics together. Thus, in the case of the section, we get two circles with the segment as a common chord (endpoints of the segment are excluded). Hereafter, we call them $\alpha$ -- isoptic circles.\end{rem}
Although the name "isoptic curve" was suggested by Taylor  in 1884 (\cite{T}), reference to former results can be found in \cite{yates}. In the obscure history of isoptic curves, we can find the names of la Hire (cycloids 1704) and Chasles (conics and epitrochoids 1837) among the contributors of the subject. A very interesting table of isoptic and orthoptic curves is introduced in \cite{yates}, unfortunately without any exact reference of its source. However, recent works are available on the topic, which shows its timeliness. In \cite{harom} and \cite{tizenketto}, the Euclidean isoptic curves of closed strictly convex curves are studied using their support function.
Papers \cite{Kurusa,Wu71-1,Wu71-2} deal with Euclidean curves having a
circle or an ellipse for an isoptic curve. Further curves appearing as isoptic curves are well studied in Euclidean plane geometry $\mathbf{E}^2$, see e.g. \cite{Loria,Wi}.
Isoptic curves of conic sections have been studied in \cite{H} and \cite{Si}. There are results for Bezier curves by Kunkli et al. as well, see \cite{Kunkli}. Many papers focus on the properties of isoptics, e.g. \cite{nyolc,MM3,het}, and the references  therein. There are some generalizations of the isoptics as well \textit{e.g.} equioptic curves in \cite{Odehnal} by Odehnal or secantopics in \cite{tizenegy, tiz} by Skrzypiec.

We can extend the very first question to the space: "What is the locus of points where a given segment subtends a given angle?" 
Or a question equivalent to the former: "For the given spatial points $A$ and $B$, what is the locus of the points $P$ for which 
the internal angle at $P$ of the triangle $ABP\triangle$ is a given angle?" We use this to define the $\alpha$ -- isoptic surface of a Euclidean spatial segment.

\begin{defn} The $\alpha$ -- isoptic surface of a Euclidean spatial segment $\overline{A_1A_2}$ is the locus of points $P$ for which the internal angle 
at $P$ in the triangle, formed by $A_1,$ $A_2$ and $P$ is $\alpha.$ If $\alpha$ is the right angle, then it is called the Thaloid of $\overline{A_1A_2}.$
\label{defisoszakaszE}
\end{defn}

It is easy to see in the Euclidean space that:
\begin{thm}
The locus of points in the Euclidean space from where a given segment subtends a given angle $\alpha$ $(0<\alpha <\pi)$ or $\pi-\alpha$ 
is a self-intersecting torus obtained by rotating the $\alpha$ -- isoptic circles drawn in any plane containing the section around the line of the section. $\square$
\end{thm}
\begin{rem} 
\begin{enumerate}
\item The torus in the above theorem contains both the isoptic surface for the given angle and the supplementary angle. In this case, we can easily separate the $\alpha$ -- and $\pi-\alpha$ -- isoptic surfaces along the self-intersection.
Specifically, the orthoptic surface is a sphere whose diameter is the section. We can call this the Thaloid of the segment.
\item There is no point in examining the isoptic surface defined in the above way for other spatial curves, because if the curve is not of constant 
0 curvature, then there is an external point from which the curve and the point cannot be fitted into a plane. In this case, the above definition needs to be generalized.
\end{enumerate}
\end{rem}
For further isoptic surfaces in Euclidean geometry, see \cite{CsSz4,CsSz5}, where we extend the definition of isoptic surfaces to other spatial objects.
The notion of isoptic curve can be extended to the other planes of constant curvature (hyperbolic plane $\mathbf{H}^2$ and spherical plane $\mathbf{H}^2$). We studied these questions in \cite{Cs-Sz14} and \cite{Cs-Sz20}. 

\section{Translation-like isoptic surfaces in $\NIL$}
In the rest of this study, we will focus on the isoptic surface of the 
translation-like segment in $\NIL$ geometry, which in the projective model is far from straight, but a parametric curve described in \ref{2.10}. 
We can make the following definition along the lines of the Definition \ref{defisoszakaszE}.

\begin{defn} The $\NIL$ translation-like $\alpha$ -- isoptic surface of a translation-like segment $\overline{A_1A_2}$ is the locus of points $P$ 
for which the internal angle at $P$ in the translation-like triangle, formed by $A_1,$ $A_2$ and $P$ is $\alpha.$ If $\alpha$ is the right angle, 
then it is called the translation-like Thaloid of $\overline{A_1A_2}.$
\label{defisoszakaszN}
\end{defn}

We emphasize here that the section itself does not appear in our calculations, 
we only deal with the endpoints. 
An interesting question beyond this study is how the ruled surface, {or more precisely in this case, how the triangular surface}  
looks like generated by the curves drawn from 
the outer point to all points 
of the section. Thus the angle can really be considered planar {in $\NIL$ sense} or any 
non desirable intersection occurs between the segment and the rays. 
The section itself and the rays can be translation-like or geodesic-like as well.
Some of these questions arise in \cite{Sz22}.

We can assume by the homogeneity of the $\NIL$ geometry that one of its endpoints coincide with the origin $A_1=E_0=(1,0,0,0)$ and the other is $A_2(1,a,b,c)$. Considering a point $P(1,x,y,z)$, we can determine the angle $A_1PA_2\angle$ along the 
procedure described in the previous section.

We apply $\bT_{P}^{-1}$ to all three points. 
This transformation preserves the angle $A_1PA_2\angle$ and pulls
back $P$ to the origin, 
hence the angle in question seems in real size. 
We get $\bT_{P}^{-1}$ by replacing $x_2,$ $y_2$ and $z_2$ in (3.2) 
with $x,$ $y$ and $z$ respectively:

\begin{equation}
\bT_{A_3}^{-1}=
\begin{pmatrix}
1 & -x & -y & xy -z \\
0 & 1 & 0 & 0 \\
0 & 0 & 1 & -x \\
0 & 0 & 0 & 1 , 
\end{pmatrix}
\end{equation}

\begin{equation}
\begin{gathered}
\bT^{-1}_{P}(P)=(1,0,0,0); \\
\bT^{-1}_{P}(A_1)=(1,-x,-y,xy -z);\\
\bT^{-1}_{P}(A_2)=(1,a-x,b-y,-bx +xy +c -z).
\label{eltoltak}
\end{gathered}
\end{equation}

According to \ref{2.10} and \ref{2.11}, the tangent of the translation curve between the origin and a point $T(1,x,y,z)$ at the origin can be obtained by the following formulas:

\begin{equation}
\bt=(u,v,w)=\left(x,y,z-\dfrac{x y}{2}\right)
\label{erintfor}
\end{equation}

Let us denote with $\bt_1$ and $\bt_2$ the tangents of the translation curves to $\bT_{P}^{-1}(A_1)$ and $\bT_{P}^{-1}(A_2)$ from the origin $E_0=\bT_{P}^{-1}(P)$ at the origin. We can calculate these tangents by applying \ref{erintfor} to \ref{eltoltak}.

\begin{equation}
\begin{gathered}
\bt_{1}=\left(-x,-y,\dfrac{x y}{2}-z\right)\\
\bt_{2}=\left(a-x,b-y,c-z+\dfrac{(x+a)(y-b)}{2}\right) 
\end{gathered}
\end{equation}

Finally, fixing the angle of the $\bt_1$ and $\bt_2$ to $\alpha,$ we get the {\it translation-like} $\alpha$ -- isoptic surface of $\overline{A_1A_2}.$ 

\begin{thm}
Given a {\it translation-like} segment in the $\NIL$ geometry by its endpoints $A_1=(1,0,0,0)$ and $A_2=(1,a,b,c).$ Then the translation-like 
$\alpha$ -- isoptic surface of the translation-like segment $\overline{A_1A_2}$ have the implicit equation:
{\footnotesize
\begin{equation}
\cos(\alpha)=\displaystyle {\frac {x\left(x-a \right) + y\left( y-b \right)+ \left(z- \frac{xy}{2} \right) \left( z-c-\dfrac{(x+a)(y-b)}{2}\right)   }
{ \sqrt{\left( {x}^{2}+{y}^{2}+ \left( z-\frac{xy}{2} \right) ^{2} \right) \left(  \left( x-a \right) ^{2}+ \left( y-b \right) ^{2}+ \left(z-c-\dfrac{(x+a)(y-b)}{2}\right)^{2} \right) 
}}}
\label{eq:niliso}
\end{equation}$\square$}
\end{thm}

\begin{figure}
\centering
\includegraphics[width=13cm]{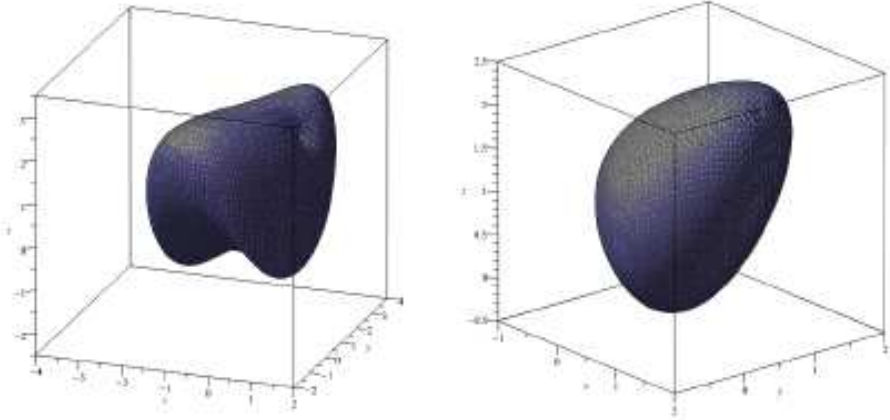}
\caption{Isoptic surface with $A_2=(1,1,1,2)$ and $\alpha=\frac{\pi}{3}$ (left) and $\frac{\pi}{2}$ (right).}
\label{fig:Fig5}
\end{figure}
On Fig.~\ref{fig:Fig5}, one can see some isoptic surfaces to a general translation-like segment in $\NIL$ geometry. The left side shows the isoptic surface to an acute angle, 
the right side shows the translation-like Thaloid of the same segment.

Let us examine the special case when the {endpoints of the segment are situated on the $z$ axis, i.e. $A_1=(1,0,0,0)$ and $A_2=(1,0,0,c).$ 
In this case, the translation-like segment looks like a Euclidean segment in the model. Replacing in \ref{eq:niliso} $a,$ $b$ and $\alpha$ with $0,$ $0$ and $\frac{\pi}{2}$ we get the following equation:
\begin{equation}
x^2+y^2+\left(z-\dfrac{xy}{2}\right)\left(z-\dfrac{xy}{2}-c\right)=0
\label{5.6}
\end{equation}
Or, after some transformation:
\begin{equation}
x^2+y^2+\left(z-\dfrac{xy}{2}-\frac{c}{2}\right)^2=\dfrac{c^2}{4}
\label{5.7}
\end{equation}

Now, applying $\bT^{-1}_F$ to all points of this equation, where $F(1,0,0,c/2),$ we obtain 

  \begin{equation}
     \begin{gathered}
       \bT^{-1}_F(P)=(1;x',y',z')=(1;x,y,z)
       \begin{pmatrix}
         1&0&0&-\frac{c}{2} \\
         0&1&0&0 \\
         0&0&1&0 \\
         0&0&0&1 \\
       \end{pmatrix}
       =\left(1;x,y,z-\frac{c}{2}\right) 
      \end{gathered} \label{5.8}
     \end{equation} 
		\begin{equation}
		(x')^2+(y')^2+\left(z'-\dfrac{xy}{2}\right)^2=\dfrac{c^2}{4}
		\label{5.9}
		\end{equation}
Comparing \ref{5.9} equation with \ref{2.13}, we can claim the following lemma:

\begin{lem}
Given a {\it translation-like} segment in the $\NIL$ geometry by its endpoints $A_1=(1,0,0,0)$ and $A_2=(1,a,b,c).$ Then the {translation-like Thaloid} 
of this line segment is a $\NIL$ sphere without the endpoints of the segment if and only if $a=b=0.$ Then the center of the $\NIL$ sphere is $F(1,0,0,c/2)$ and its radius is $\dfrac{|c|}{2}.$
\label{lemiso}
\end{lem}
\textbf{Proof:} To prove Lemma \ref{lemiso}, we need further consideration to the other direction, not covered by the calculations above. 
Due to lengthy calculations and formulas, we only present here the outline of the proof.
First, we need to find the midpoint $F$ of the translation section $\overline{A_1A_2}$ and its length, half of which will be the radius 
of the sphere. Then, using $\bT^{-1}_F$ and \ref{2.13}, we write the equation of the translation sphere whose diameter is $\overline{A_1A_2}.$ 
We compare this to the numerator of the right side in equation \ref{eq:niliso} (substituting $\alpha=\pi/2$). Finally, considering the 
difference of the two equations, we get that it will be $0$ if and only if $ay-bx=0,$ or $ay-bx=2ab-4c$ which is equivalent to $a=b=0.$ 
$\square$

On Fig.~\ref{fig:Fig6} we can see the {translation-like Thaloid} related to the translation segment $A_1=(1,0,0,0)$ and $A_2=(1,0,0,4),$ as well as a sphere with center $F=(1,0,0,2)$ and radius 2. 
\begin{figure}
\centering
\includegraphics[width=6cm]{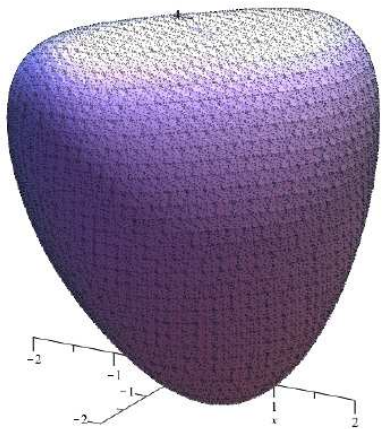}
\caption{Isoptic surface (translation-like Thaloid) of translation segment $\overline{A_1A_2}$ where $A_1=(1,0,0,0)$, $A_2=(1,0,0,4)$ and $\alpha=\frac{\pi}{2}.$}
\label{fig:Fig6}
\end{figure}


%
\end{document}